\documentclass{amsart}
\usepackage{amssymb}
\usepackage{graphicx}
\def\({\left(}
\def\){\right)}

\newtheorem{lema}{Lemma}[section]

\newtheorem*{teorema*}{Theorem}

\newtheorem{remark}[lema]{Remark}

\newtheorem{lemma}{Lemma}[section]
\newtheorem{corollary}[lema]{Corollary}
\newtheorem{theorem}[lema]{Theorem}

\newtheorem{proposition}[lema]{Proposition}

\newtheorem{definition}[lema]{Definition}

  {\par \hfill \fbox{}}
  {\par \hfill \fbox{}}
  {\par \hfill }
  {\par \hfill }

\def\beq{\begin{equation}}
\def\eeq{\end{equation}}
\def\beginpf{\noindent{\bf Proof.} \quad}
\def\endpf{\qquad\hfill\rule{2.2mm}{2.2mm}\break}

\def\epsilon{\varepsilon}

\begin{document}

\title[quaternionic Fock space]{ quaternionic Fock space  on slice  hyperholomorphic functions }
\author{Sanjay Kumar}
\address{Department of Mathematics, Central University of Jammu,
Jammu 180 011, INDIA.} \email{sanjaykmath@gmail.com}
\author{S. D. Sharma}
\address{Department of Mathematics, Central University of Jammu,
Jammu 180 011, INDIA.} \email{somdatt\_jammu@yahoo.co.in}
\author{Khalid Manzoor}
\address{Department of Mathematics, Central University of Jammu,
Jammu 180 011, INDIA.}
\email{khalidcuj14@gmail.com}

%\thanks{The second author is supported by NBHM (DAE) Post-Doctoral Fellowship grant no. 2/40(32)/2009-R\&D-II/1337.\\
%The third author is supported by the National Natural Science Foundation of China  grant nos. 10971153, 10671141}

\date{\today.}
%\dedicatory{}
%\commby{}

\begin{abstract}
In this paper,  we define the quaternionic  Fock spaces $\mathfrak{F}_{\alpha}^p$ of entire  slice hyperholomorphic functions in a quaternionic unit ball $\mathbb{B}$ in $\mathbb{H}.$ We also study  growth  estimate  and various results  of entire slice regular functions in these  spaces. The work of this paper is motivated by the recent work of \cite{alpa12} and \cite{marco}.
\end{abstract}

% ----------------------------------------------------------------------
%\begin{center}
\maketitle
%\end{center}
% ----------------------------------------------------------------------

\section{Introduction}
The notation of slice hyperholomorphicity is introduced in 2006 in   \cite{gent06}  and till then a lot of works  have been done in this direction. Several function spaces like Hardy spaces,  Bergman spaces, Bloch,  Besov and Dirichlet  spaces, Pontryagin De Branges Rovnyak spaces, etc  are studied in the  slice hyperholomorphic settings, see \cite{  alpa12, alpa13,  alpa14, arco14, arco15,  colo12,   colo13,  colo131, sarf13,  marco}.    We refer to a recent survey \cite{colo15}  and the book  \cite{colo11} for details  information and references for the systematic development of slice hyperholomorphic functions and their  applications.  The Fock spaces in the slice hyperholomorphic settings were studied by D. Alpay, F. Colombo, I. Sabadini,  \cite{{alpa14}}. The Fock spaces are fundamental for their role in quantum mechanics, see \cite{adle95, arco15,  zhu12} and references therein.
  By symbol  $\mathbb{H} =\{x_0+x_1i+x_2j+x_3k :x_l \in \mathbb{R} ~\mbox{for}~0\leq l \leq 3\},$ we  denote the  set of    4-dimensional non-commutative real algebra of quaternions    generated by    imaginary units
$i,j,k$  such that $i^2=j^2=k^2=-1$ and $ ij=-ji=k, jk=-kj=i,  ki=-ik=j.$  The Euclidean norm of a quaternion $q$ is given by     $|q| =\sqrt{q\bar q}=\sqrt{\bar q q}=\sqrt{\sum_{l=0}^3}x_l^2,~ \mbox{for}~ x_l \in \mathbb{R},$  where $\bar q=Rel(q)-Im(q)= x_0-(x_1i+x_2j+x_3k),$ denote  the congugate of $q.$   The multiplicative inverse of non-zero quarternion $q$ is given by $\displaystyle\frac{\bar q}{|q|^2}.$  The set $\mathbb{S}=\{q\in \mathbb{H}:  q=x_1i+x_2j+x_3k ~\mbox{and}~x_1^2+x_2^2+x_3^2=1\}$  represents the two-dimensional  unit sphere of purely imaginary quaternions. Any element $I\in \mathbb{S}$ is such that  $I^2=-1.$ This implies that the elements of S are  imaginary units.  The quaternion is considered as the union of  complex  plane  $\mathbb{C}_I=\mathbb{R}+\mathbb{R}I$  (also called slices)  each one is identified by an imaginary unit $I\in \mathbb{S}.$ Let   $\Omega_I= \Omega \cap \mathbb{C}_I,$ for  some 
domain $\Omega$ of $\mathbb{H}.$
 For any  quaternion $q$ we can write $q=x_0+x_1i+x_2j+x_3k=x_0+Im(q)=x_0+|Im(q)|I_q=x+yI_q,$
 with
$ I_q=\displaystyle\frac{Im(q)}{|Im(q)|}$ if  $|Im(q)|\neq 0,$ otherwise we take arbitrary $I$ in $\mathbb{S}.$

%%%%%%%%%%%%%%%%%%%%%%%%%%%%%%%%%%%%%%%
Here, we begin with some  basic results  in the quaternionic-valued slice regular functions.
\begin{definition}
Let $\Omega$ be a domain in $\mathbb{H}$.  A real differentiable  function $f:\Omega \to \mathbb{H}$ is said to be the (left) slice regular or slice hyperholomorphic if for any $I\in \mathbb{S}, f_I$  is holomorphic in $\Omega_I$, i.e:,  $$ \left(\frac{\partial }{\partial x_0}+I\frac{\partial }{\partial y}\right)f_I(x_0+yI)=0,$$  where $f_I$ denote the restriction of $f$ to 
  $\Omega_I.$ The  class of slice regular function on $\Omega$ is denoted  by $SR(\Omega).$
\end{definition}

%%%%%%%%%%%%%%%%%%%%%%%%%%%%%%%%%%%%%%%%%
For slice regular functions, we have the following useful result.
\begin{theorem} \cite[Theorem  2.7]{gent07}
A function $f: \mathbb{B} \to \mathbb{H}$ is said to be slice regular if and only if it has  a power series of the form 
\beq\label{eq:19}
f(q)=\displaystyle\sum_{n=0}^\infty q^na_n, ~\mbox{where}~a_n=\frac{1}{n!}\frac{\partial^n f(0)}{\partial x^n}
\eeq
 converging uniformly on $\mathbb{B}.$
\end{theorem}
%%%%%%%%%%%%%%%%%%%%%%%%%%%%%%

\begin{lemma}\label{eq:1}
 \cite[Lemma 4.1.7 ]{colo11}(Splitting Lemma)  If $f$ is a slice regular function on the domain $\Omega,$ then for any 
$i,j\in  \mathbb{S},$  with $i\bot j$ there exists two holomorphic functions $F,L:\Omega_I \to\mathbb{C}_I$ 
such that 
\beq \label{eq:2}
f_I(z)=F(z)+L(z)J ~\mbox{for any }~z=x+yI.
\eeq
\end{lemma}
%%%%%%%%%%%%%%%%%%%%%%%%%%%%%%%%%%
%%%%%%%%%%%%%%%%%%%%%%%%%%%%%%%%%%%%%%%%%%
\noindent One of the most important property of the slice regular functions is their Representation Formula.  It only  holds on the open sets which are stated below.
\begin{definition}
 Let $\Omega$ be an open set in  $\mathbb{H}$. We say  $\Omega$ is axially symmetric if for any  $q=x+yI_q \in \Omega$ all the elements $x+yI$ are contained in $\Omega$, for all $I\in  \mathbb{S}$ and $\Omega$ is said to be slice domain if $\Omega\cap \mathbb{R}$ is non empty and  $\Omega\cap \mathbb{C}_I$ is a domain in $\mathbb{C}_I$ for all $I\in \mathbb{S}.$
\end{definition}
\begin{theorem} \label{eq:115}
  \cite[Theorem  4.3.2 ]{colo11} (Representation Formula)  Let  $f$ be a  slice regular function in the domain $\Omega\subset  \mathbb{H}.$ Then for any $j\in  \mathbb{S}$ and for all $z=x+yI \in \Omega, $   $$f(x+yI)=\frac{1}{2}\{(1+IJ)f(x-yI)+(1-IJ)f(x+yI)\}.$$
\end{theorem}
%%%%%%%%%%%%%%%%%%%%%%%%%%%%%%%%%%%%%%%%%%%%%%%%%%%%
\begin{remark}
Let $I,J$ be orthogonal imaginary units in $\mathbb{S}$ and $\Omega$ be  an axially  symmetric slice domain. Then the Splitting Lemma and the  Representation formula  generate a class of  operators on the slice regular functions as follows:

$$Q_I:SR(\Omega)\to hol({\Omega}_I)+hol({\Omega}_I)J$$ $$Q_I: f\mapsto  f_1+f_2J$$

$$P_I: hol({\Omega}_I)+hol({\Omega}_I)J \to SR(\Omega)$$
$$PI[f](q)=P_I[f](x+yI_q)=\frac{1}{2}[(1-I_qI)f(x+yI)+(1+I_qI)f(x-yI)].$$ Also, $$P_I\circ Q_I=\emph{I}_{SR(\Omega)}~ \mbox{and}~Q_I\circ P_I=\emph{I}_{SR( hol({\Omega}_I)+hol({\Omega}_I))},$$ where $\emph{I}$ is an identity operator.

\end{remark}

Since pointwise product of functions does not preserve slice regularity (see \cite{colo11}) a new multiplication operation for regular functions is defined. In the special case of power series, the regular product (or $\star-$product) of $f(q) = \sum_{n= 0}^{\infty}q^{n}a_{n}$ and  $g(q) = \sum_{n= 0}^{\infty}q^{n}b_{n}$ is $$ f \star g(q) =  \sum_{n \ge 0} q^{n} \sum_{k = 0}^{n} a_{k}b_{n-k}. $$
The $\star-$product is related to the standard pointwise product by the following formula.
\begin{theorem} 
  \cite[Proposition   2.4 ]{arco14}
Let $f, g$ be regular functions on $\mathbb{B}.$ Then
$f \star g(q) = 0 $ if $f(q)  = 0 $ and $ f(q) g(f(q)^{-1} qf(q))$ if $ f(q) \neq 0.$ 
The reciprocal $f^{-\star}$ of a regular function $f (q) =\sum_{n= 0}^{\infty}q^{n}a_{n}$
with respect to the $\star-$product is
$$f^{\star}(q) = \frac{1}{f \star f^{c}(q)} f^{c}(q), $$ where $f^{c}(q) =\sum_{n= 0}^{\infty}q^{n}\overline{a_{n}} $ is the regular conjugate of $f.$ The function $f^{-\star}$
is regular on $\mathbb{B}\setminus(q \in \mathbb{B} | f \star f^{c} (q) = 0) $ and $f  \star f^{-\star} = 1$ there.
\end{theorem}
%%%%%&&&&&&&&&&&&&&&&&&&&&&&&&&&&&&&&&&&&&&&&&&&&&&&&&&&&&&&&&&&&&&&&&&&&&&&&&&&&&&&
\section{Fock spaces}
\label{sec:2}

 In  this section, we study some basic properties of Fock spaces in the slice hyperholomorphic  settings. Fock spaces of holomorphic functions are discussed in  details in the book \cite{zhu12}. The slice hyperholomorphic quaternionic Fock spaces are studied in \cite{alpa14}. Let $dA$ be the normalized area measure on  $\mathbb{C}.$
 For $0<p< \infty,$  the   Fock space  $\mathfrak{F}_{p, \mathbb{C}}$ is defined as the space of entire functions  $f: \mathbb{C} \to  \mathbb{C}$  such that $$\displaystyle\frac{\alpha p}{2\pi}\int_{\mathbb{C}}\left|f(z)e^{\frac{-\alpha}{2}|z|^2}\right|^p dA(z) <\infty,$$
where  $ z  \in \mathbb{C}$ and  $dA(z) =\frac{1}{\pi} dxdy,~z=x+jy,~x, y \in \mathbb{R}.$
 Let $\mathbb{B}(0,1)=\mathbb{B}=\{q=x+yI_q: |q|<1\}$ be the quaternionic  unit ball  centered at origin in $\mathbb{H}$ and $\mathbb{B}\cap \mathbb{C}_I=\mathbb{B}_I$ denote   unit disk in the complex plane   $\mathbb{C}_I$ for $ I\in \mathbb{S}.$     A function slice regular on the quaternionic space  $\mathbb{H}$  is called  all  slice regular and have power series representation of the form (\ref{eq:19}) converging everywhere in  $\mathbb{H}$  and uniformily on the compact subsets of  $\mathbb{H}$. Let  $SR(\mathbb{H})$ denote  the space of  entire slice regular functions on the unit ball $\mathbb{B}.$
  Here we begin with the following definitions.
\begin{definition}
 For $0<p<\infty$ and  $I\in \mathbb{S},$  the quaternionic right linear space of entire slice regular functions $f$ is said to be the quaternionic slice  regular Fock  space on the unit ball $\mathbb{B},$ if for any $q\in \mathbb{B}$
$$\displaystyle\frac{\alpha p}{2\pi}\sup_{I\in \mathbb{S}}\int_{\mathbb{B}_I}\left|f(q)e^{\frac{-\alpha}{2}|q|^2}\right|^p dA_I(q) <\infty,$$that is,
$$\mathfrak{F}_{\alpha}^p=\{f\in SR(\mathbb{H}):\frac{\alpha p}{2\pi}\sup_{I\in \mathbb{S}}\int_{\mathbb{B}_I}\left|f(q)e^{\frac{-\alpha}{2}|q|^2}\right|^p dA_I(q) <\infty\},$$ where $ dA_I(q)$ denote the normalized differential area in the complex plane $\mathbb{C}_I$ such that area of $\mathbb{B}_I$ is equal to one
 and is M\"{o}bius invariant measure on $\mathbb{B}$ with   norm given  by  $$\|f\|_{\mathfrak{F}_{\alpha}^p}=\left( \frac{\alpha p}{2\pi}\sup_{I\in \mathbb{S}}\int_{\mathbb{B}_I}\left|f(q)e^{\frac{-\alpha}{2}|q|^2}\right|^p dA_I(q): ~q=x+yI_q \in \mathbb{B}\right)^{\frac{1}{p}} .$$
\end{definition}
 \noindent By $\mathfrak{F}_{\alpha, I}^p,$ we denote the quaternionic right linear space of entire slice regular functions on  $\mathbb{B}$ such that  $$\displaystyle\frac{\alpha p}{2\pi}\int_{\mathbb{B}_I}\left|f(z)e^{\frac{-\alpha}{2}|z|^2}\right|^p dA_I(z) <\infty.$$ Furthermore, for each function $f\in \mathfrak{F}_{\alpha,I}^p,$ we define $$\|f\|_{\mathfrak{F}_{\alpha,I}^p} = \left(\frac{\alpha p}{2\pi}\int_{\mathbb{B}_I}\left|f(z)e^{\frac{-\alpha}{2}|z|^2}\right|^p dA_I(z)~:z=x+yI\in  \mathbb{B}\cap\mathbb{C}_I\right)^{\frac{1}{p}}.$$

\begin{remark} \cite[P. 499 ]{marco}\label{eq:100}
Let $I\in \mathbb{S}$ be such that $J\bot I.$ Then there exist holomorphic functions  $f_1, f_2: \mathbb{B}_I \to \mathbb{C}_I$ such that $f_I=Q_I[f]=f_1+f_2J $ for some holomorphic map $Q_I[f]$ in complex variable $z \in \mathbb{B}_I$. Then by  \cite[Remark 3.3]{marco}, we have 
$$\begin{array}{ccl}
\left|f_l(z)e^{\frac{-\alpha}{2}|z|^2}\right|^p&\leq &\left|f(z)e^{\frac{-\alpha}{2}|z|^2}\right|^p \\
&\leq& 2^{max\{0,p-1\}}\left|f_1(z)e^{\frac{-\alpha}{2}|z|^2}\right|^p+2^{max\{0,p-1\}}\left|f_2(z)e^{\frac{-\alpha}{2}|z|^2}\right|^p.
\end{array}$$
The condition  $f\in \mathfrak{F}_{\alpha,I}^p$  is equivalent to $f_1$ and $f_2$ belonging to  one dimensional   complex Fock space.
\end{remark}
%%%%%%%%%%%%%%%%%%%%%%%%%%%%%%%%%%%%%%%%%%%%%
\begin{proposition}\label{eq:160}
Suppose   $I\in  \mathbb{S}$ and  $\alpha>0.$ Then $f\in  \mathfrak{F}_{\alpha,I}^p, \;  p>1$ if and only if $f\in  \mathfrak{F}_{\alpha}^p.$ Moreover, the spaces $( \mathfrak{F}_{\alpha,I}^p,\|.\|_{ \mathfrak{F}_{\alpha,I}^p})$ and $( \mathfrak{F}_{\alpha}^p,\|.\|_{ \mathfrak{F}_{\alpha}^p})$ have equivalent norms. More precisely, one has  $$\|f\|_{\mathfrak{F}_{\alpha,I}^p}^p\leq \|f\|_{\mathfrak{F}_{\alpha}^p}^p\leq 2^p\|f\|_{\mathfrak{F}_{\alpha,I}^p}^p.$$
\end{proposition}
\beginpf Let  $f\in  \mathfrak{F}_{\alpha}^p.$ Since $\mathbb{B}_I\subset \mathbb{B}.$ Then by definition, $\|f\|_{\mathfrak{F}_{\alpha,I}^p}^p\leq \|f\|_{\mathfrak{F}_{\alpha}^p}^p$ which implies $\mathfrak{F}_{\alpha}^p\subset \mathfrak{F}_{\alpha,I}^p.$ Now, let  $f\in  \mathfrak{F}_{\alpha,I}^p.$  For $q=x+yI_q \in \mathbb{B}$ with $I_q=\frac{Im(q)}{|Im(q)|}$ and $z=x+yI \in \mathbb{B}_I $ and as  $|q|=|z|.$ Then thanks to Representation Formula for slice regular  functions, we have
$$\begin{array}{ccl}
\displaystyle \frac{\alpha p}{2\pi}\int_{\mathbb{B}_I}\left|f(q)e^{\frac{-\alpha}{2}|q|^2}\right|^pdA_I(q)&=&\displaystyle \frac{\alpha p}{2\pi}\int_{\mathbb{B}_I}\frac{1}{2}|(1-I_qI)(f(z)e^{\frac{-\alpha}{2}|z|^2})\\
&+&\displaystyle (1+I_qI)(f(\bar z)e^{\frac{-\alpha}{2}|\bar z|^2})|^pdA_I(z)\\
&\leq& 2^{max\{0, p-1\}}\displaystyle \frac{\alpha p}{2\pi} \int_{\mathbb{B}_I}\left|f(z)e^{\frac{-\alpha}{2}|z|^2}\right|^pdA_I(z)\\
&+& \displaystyle2^{max\{0,p-1\}}\displaystyle \frac{\alpha p}{2\pi} \int_{\mathbb{B}_I}\left|f(\bar z)e^{\frac{-\alpha}{2}|\bar z|^2}\right|^pdA_I(\bar z).
\end{array}$$
Hence on taking supremum over all $I\in \mathbb{S},$ we have 
$$\begin{array}{ccl}
 \|f\|_{\mathfrak{F}_{\alpha}^p}^p&\leq& 2^{max\{0, p-1\}}\displaystyle \frac{\alpha p}{2\pi}\left( \int_{\mathbb{B}_I}\left|f(z)e^{\frac{-\alpha}{2}|z|^2}\right|^pdA_I(z)+ \int_{\mathbb{B}_I}\left|f(\bar z)e^{\frac{-\alpha}{2}|\bar z|^2}\right|^pdA_I(\bar z) \right)\\
&\leq & 2^{p-1} 2 \|f\|_{\mathfrak{F}_{\alpha,I}^p}.
\end{array}$$
\endpf
%%%%%%%%%%%%%%%$$$$$$$$$$$$$$$$$$$$$$$$$$$$
 We can easily   prove the following results.
\begin{proposition}\label{eq:153}
  Suppose  $p>1,\alpha >0.$ If  $f \in SR(\mathbb{H}),$ then  following statements are equivalent:\\
(a) $f\in \mathfrak{F}_{\alpha}^p;$\\
(b) $f\in \mathfrak{F}_{\alpha,I}^p $
 for some $I\in \mathbb{S}.$
\end{proposition}
%%%%%%%%%%@@@@@@@@@@@@@@@@@@@@@@@@@@@@@@@@@@@@@@@@@@@@@@
\noindent
\begin{proposition}
Let   $I,J\in \mathbb{S},  p>1$ and $\alpha >0.$ If $f\in SR(\mathbb{H}),$ then  $f\in \mathfrak{F}_{\alpha, I}^p$  if and only if   $f\in \mathfrak{F}_{\alpha, J}^p.$
\end{proposition}

\begin{proposition}\label{eq:129}
The space $\mathfrak{F}_{\alpha}^p,~ p>1$ and $\alpha>0$  is complete.
\end{proposition}
\beginpf
Let $\{f_m\}_{m\in \mathbb{N}}$ be a Cauchy sequence in $\mathfrak{F}_{\alpha}^p.$ Then, for $I\in \mathbb{S},$  $\{f_m\}$ is Cauchy sequence in  $\mathfrak{F}_{\alpha,I}^p.$ Let $J\in \mathbb{S} $  be such that $J\bot I$ and let $f_{m,1}, f_{m,2}$ be holomorphic functions such that $f_I=Q_I[f]=f_{m,1}+f_{m,2}J.$ Since $\{f_{m,1}\}_{m\geq 0}$ and  $\{f_{m,2}\}_{m\geq 0}$ are  Cauchy sequences in the complex Fock space $\mathfrak{F}_{\alpha, \mathbb{C}_I}^p$ and  the fact that $\mathfrak{F}_{\alpha, \mathbb{C}_I}^p$ is complete, so  we conclude that, there exist functions $f_l\in \mathfrak{F}_{\alpha, \mathbb{C}_I}^p$ such that  each $f_{m,l}\to f_l $ as $m\to \infty$ for $l=1,2.$ Now set  $f=P_I(f_1+f_2J).$ Therefore, $$\|f_m-f\|^p_{\mathfrak{F}_{\alpha,I}^p}\leq \|f_{m,1}-f_1\|^p_{\mathfrak{F}_{\alpha,\mathbb{C}_I}^p}+\|f_{m,2}-f_2\|^p_{\mathfrak{F}_{\alpha,\mathbb{C}_I}^p}\to 0~~\mbox{as}~~m\to \infty.$$ This implies that  $f_m \to f$ in $\mathfrak{F}_{\alpha,I}^p.$  Hence $f\in \mathfrak{F}_{\alpha,I}^p$ and so $f\in \mathfrak{F}_{\alpha}^p$. Thus, the    slice regular Fock  space $\mathfrak{F}_{\alpha}^p$ is complete.
\endpf
%%%%%%%
\begin{remark}
  If we write  $$d\lambda_{\alpha, I}(q)=\displaystyle\frac{\alpha}{\pi}e^{-\alpha |q|^2}dA_I(q);  q=x+yI_q \in \mathbb{B},$$ then the slice regular Fock space  has the structure of  quaternionic  Hilbert space  with  their inner product  $\langle . , . \rangle_{\alpha}$ defined by $$\langle f , g\rangle_{\alpha}=\displaystyle \int_{\mathbb{B}_I}f(q)  \overline {g(q)}d\lambda_{\alpha, I}(q)$$ for $f,g\in \mathfrak{F}_{\alpha}^p.$
\end{remark}
\begin{proposition}\label{eq:130}
  On the  slice regular Fock  $\mathfrak{F}_{\alpha}^p,$ the function    $\langle . , . \rangle_{\alpha}$ is a   quaternionic  right linear  inner product,  i.e., for all $f,g,h\in \mathfrak{F}_{\alpha}^p$ and $a \in \mathbb{H}$
\begin{itemize}
\item[$(i)$] positivity: $\langle f , f\rangle_\alpha\geq 0$ and $\langle f , f\rangle_\alpha = 0$ if and only if $f=0;$ 
\item[$(ii)$]  quaternionic hermiticity: $\langle f , g\rangle_\alpha=\overline{\langle  g, f\rangle_\alpha};$
\item[$(iii)$] right linearity:  $\langle f,g a+h\rangle_\alpha=\langle f , g\rangle_\alpha a+\langle f , h\rangle_\alpha.$
\end{itemize}
\end{proposition}
%%%%%%%%$$$$$$$$$$$$$$$$$$$$$$$$$$$$$$$$$$$$$$$$$$$$$$$$$$$$$$$$$$$$$$$$$$
\begin{proposition}
For $ p>1$ and $\alpha>0,$  The space $(\mathfrak{F}_{\alpha}^p, \langle . , .\rangle_{\alpha})$ is   quaternionic  Hilbert  space.
\end{proposition}
\beginpf
 From Proposition  \ref{eq:130}, it follows that the function $\langle .  , . \rangle_{\alpha}$ is  a quaternionic  right linear  inner product and   Proposition  \ref{eq:129}  shows that the slice regular Fock space  is complete. 
\endpf
%%%%%%%%%%%%%%%%%%%%%%%%%%%%%%%%%%%%%%%%%%%%%%%%
%$$$$$$$$$$$$$$$$$$$$$$$$$$$$$$$$$$$$$$$$$$$$$$$$$$$$$$$$$$
%%%%%%&&&&&&&&&&&&&&&*************************************************************************
%%%%%%%%###################################################################
\begin{remark} By $L^{p}(\mathbb{B}_I, d\lambda_{\alpha, I}, \mathbb{H}),$  we define   the set of functions $g: \mathbb{B}_I \to \mathbb{H} $ such that $$\displaystyle\int_{\mathbb{B}_I} |g(w)|^p d\lambda_{\alpha, I}(w)<\infty,$$ where $d\lambda_{\alpha, I}(w)=\frac{\alpha}{\pi}e^{-\alpha|z|^2}dA_I(w)$ for $\alpha>0$  is called  the Gaussian probability measure. Note that for $J\in \mathbb{S}$ with $J\bot I$ and $g=g_1+g_2J$ with $g_1,g_2:\mathbb{B}_I\to \mathbb{C}_I,$ then $g\in L^{p}(\mathbb{B}_I, d\lambda_{\alpha, I}, \mathbb{H})$ if and only if  $g_1 ,g_2 \in L^{p}(\mathbb{B}_I, d\lambda_{\alpha, I}, \mathbb{C}_I).$\\ Clearly, $\mathfrak{F}_{\alpha}^p$ is a closed subspace of  $L^p(\mathbb{B}_I, d\lambda_{\alpha, I}, \mathbb{H}).$
In complex analysis, the reproducing kernel of complex Fock space for $p=2$ is given by $$K_{\alpha}^{\mathbb{C}_I}(z,w)=e^{\alpha\langle z, w\rangle};~~z,w\in \mathbb{C}_I.$$
\end{remark}
 This gives the motivation  for the following definition.
%%%%%%&&&&&&&&&&&&&&&&&&&&&&&&&&&&&&&&&&&&&&&&&&&&&&&&&&&&&&&&&&&&&&&&&&&&&&

\begin{definition}
For any $q\in \mathbb{B},$ the slice regular exponential function is given by $$e^q= \sum_{n=0}^{\infty} \frac{q^n}{n!}.$$
Let  $\displaystyle e^{zw}= \sum_{n=0}^{\infty} \frac{z^n w^n}{n}$ be a holomorphic function in variable $z$ in the complex plane $\mathbb{C}_I.$  Clearly, $e^{zw}$ is not slice regular in both  variable.  Setting $\displaystyle e^{qw}_{\star}= \sum_{n=0}^{\infty} \frac{q^n w^n}{n!},$ then we see that  the function  is left slice regular in q and right slice regular in w, where $\star$  denote  the slice regular  product. By Representation Formula, we can obtain the extension  of  function $e^{zw}$ to $\mathbb{H},$  as $$ext(e^{zw})=\frac{1}{2}\{(1-IJ)e^{zw}+(1+IJ)e^{\bar{z}w}\}=e^{qw},$$ where $q\in \mathbb{B}$ and for some arbitrary $w.$
For $I\in \mathbb{S}$ and $\alpha>0,$ we define $$B_{\alpha}(q,w)=e^{\alpha q\bar w}_{\star}~ for~ each ~q\in \mathbb{B}$$ and  is called slice regular reproducing kernel  of quaternionic  Fock space. 
\end{definition}

%%%%^^^^^^^^^^^^^^^^^^^^^^^^^^^^^^^^^^^^^^^^^^^^^^^^^^
%%%%################################################################
\begin{proposition}
For any positive integer $m,$ the set of the form $e_m(q)=q^m\sqrt{\frac{\alpha}{m}}$   is orthonormal in the  quaternionic  Fock space  $\mathfrak{F}_{\alpha}^2(\mathbb{B}).$
\end{proposition}
\beginpf
By Lemma \ref{eq:1},  we can write $f_I=f_1+f_2J$ for some $\mathbb{C}_I$-valued holomorphic functions $f_1,f_2.$ Now for any $m>0,$   we have $$\begin{array}{ccl}
\langle f, e_m\rangle_\alpha&=& \langle f_1+f_2J, e_m\rangle_\alpha=\langle f_1,e_m\rangle_\alpha+\langle f_2J,e_m\rangle_\alpha \\
&=&\displaystyle\int_{\mathbb{B}_I}f_1(z)e_m d\lambda_{\alpha, I}(z)+\displaystyle\int_{\mathbb{B}_I}f_2(z)e_m d\lambda_{\alpha, I}(z)J.
\end{array}$$
In complex plane every power series of the form $f_l(z)=\displaystyle\sum_{k=0}^\infty z^ka_{l,k},~ l=1,2$ converges uniformly on $|z|<R,$ for each $z\in \mathbb{B}_I.$\\ Therefore, we obtain
 $$\begin{array}{ccl}
<f, e_m>_\alpha&=&\displaystyle\int_{|z|<R}\displaystyle\sum_{k=0}^\infty z^ka_{1,k}e_m(z)d\lambda_{\alpha, I}(z)+\displaystyle\int_{|z|<R}\displaystyle\sum_{k=0}^\infty z^ka_{2,k}e_m(z)d\lambda_{\alpha, I}(z)J\\
&=&\displaystyle\sum_{k=0}^\infty a_{1,k} \displaystyle\int_{|z|<R}z^ke_m(z)d\lambda_{\alpha, I}(z)+\displaystyle\sum_{k=0}^\infty a_{2,k} \displaystyle\int_{|z|<R}z^ke_m(z)d\lambda_{\alpha, I}(z)J\\
&=&\displaystyle\lim_{R\to \infty} (a_{1,m}+a_{2,m}J) \displaystyle\int_{|z|<R}z^ke_md\lambda_{\alpha, I}(z)\\
&=&\displaystyle\lim_{R\to \infty}d_m\int_{\mathbb{B}_I} q^ke_m(q)d\lambda_{\alpha, I}(q),
\end{array}$$
where $d_m= a_{1,m}+a_{2,m}J.$ But in complex Fock space $\mathfrak{F}_\alpha^2(\mathbb{B}_I),$ each $a_{l,k}=0$ for $ l=1,2$ implies $d_m=0$ and so $f=0.$ Thus the sequence $\{e_m\}_{m>0}$ is complete in  $\mathfrak{F}_\alpha^2(\mathbb{B}).$
\endpf
%%%%%%%%%&&&&&&&&&&&&&&&&&&&&&&&&&&&&&&&
%%%%%%%%^^^^^^^^##############################################################
\begin{proposition}
 For some $I \in \mathbb{S},$ the slice regular  orthogonal projection on  $\mathbb{B}$ is defined  by    $T_{\alpha,I}: L^2(\mathbb{B}_I, d\lambda_{\alpha, I}, \mathbb{H}) \to  \mathfrak{F}_{\alpha}^2.$   Then  for all $q,w\in \mathbb{B},$ the integral representaion  for $T_{\alpha,I}$ is given by $$T_{\alpha,I}f(q)=\displaystyle\frac{\alpha p}{2\pi}\int_{\mathbb{B}_I}  f(w)B_{\alpha}(q,w) e^{-\alpha |w|^2}d A_I(w)$$ for all $f\in  L_2(\mathbb{B}_I, d\lambda_{\alpha, I}, \mathbb{H}), $   where $B_{\alpha}(q,w)=e^{\alpha\langle q,w\rangle}_{\star}=e^{\alpha q \bar w}_{\star}$ is reproducing kernel  for $\mathfrak{F}_{\alpha}^2.$
\end{proposition}
\beginpf
Given   $f\in  L^2(\mathbb{B}_I, d\lambda_{\alpha, I}, \mathbb{H}), $ let  $Q_I[f]=f_I$  be its restriction.  Then we write  $Q_I[f]=f_1+f_2J,$  where $J$ is an element of  $\mathbb{S}$ such that $J\bot I$ and $f_1, f_2$ are complex valued holomorphic functions. Further, if  for all $z, w\in \mathbb{B}_I,$ then  the two functions $K^{\mathbb{C}_I}_{\alpha}(z,w)$ and $B_{\alpha}(z,w)$ coincide and  from the fact that $f=\langle f,B_\alpha\rangle$ (see \cite[Theorem 3.10]{alpa14}), one conclude
$$\begin{array}{ccl}
T_{\alpha,I}f&=&\langle T_{\alpha,I}f, B_{\alpha}(.,.) \rangle_{\alpha}\\
&=&\langle T_{\alpha,I}(f_1+f_2J), K^{\mathbb{C}_I}_\alpha\rangle_{\alpha}\\
&=&\langle T_{\alpha,I}f_1, K^{\mathbb{C}_I}_\alpha\rangle_{\alpha}+\langle T_{\alpha,I}f_2J, K^{\mathbb{C}_I}_\alpha\rangle_{\alpha}\\
&=&\langle f_1, K^{\mathbb{C}_I}_\alpha\rangle_{\alpha}+\langle f_2J, K^{\mathbb{C}_I}_\alpha\rangle_{\alpha}\\
&=&\langle f_1+f_2J, K^{\mathbb{C}_I}_\alpha\rangle_{\alpha}\\
&=&\langle f,B_\alpha\rangle_{\alpha}\\
&=& \displaystyle\int_{\mathbb{B}_I}  f(w)B_{\alpha}(q,w) d\lambda_{\alpha, I}(w)\\
&=&\displaystyle\frac{\alpha p}{2\pi}\int_{\mathbb{B}_I}  f(w)B_{\alpha}(q,w) e^{-\alpha |w|^2}d A_I(w).\\
\end{array}$$
\endpf
In the next result, we give the growth rate estimation   for  entire slice regular  functions in quaternionic  Fock space.
\begin{theorem} 
Let $1< p\leq \infty$ and $\alpha>0.$  Then for every $f\in \mathfrak{F}_{\alpha}^p,$ $$\displaystyle\sup_{q\in \mathbb{B}}\left\{|f(q)|:\|f\|_{\mathfrak{F}_{\alpha}^p}\le 1\right\}\leq2 e^{\frac{\alpha}{2}|q|^2},~where~q=x+yI_q~and  ~I_q=\frac{Im(q)}{|Im(q)|}.$$
\end{theorem}
\beginpf
Let $I,J$ be the orthogonal imaginary units in two dimensional sphere $\mathbb{S}$. If  $f\in \mathfrak{F}_{\alpha}^{p},$ then by Proposition  \ref{eq:153}, $f\in \mathfrak{F}_{\alpha,I}^p$ and $\|f\|_{\mathfrak{F}_{\alpha,I}^p}\le 1.$ Now,  we can find   two holomorphic functions $f_1,f_2$ in $\mathfrak{F}_{\alpha,\mathbb{C}_I}^p$ such that $Q_I[f]=f_1+f_2J.$ By using \cite[Theorem 2.7]{zhu12}, each  $f_l$  satisfies  $\displaystyle\sup_{z\in \mathbb{B}_I}\left\{|f_l(z)|:\
\|f_l\|_{\mathfrak{F}_{\alpha,\mathbb{C}_I}^p}\le 1\right\}= e^{\frac{\alpha}{2}|z|^2}.$ Furthermore, $$\displaystyle\sup_{z\in \mathbb{B}_I}\left\{|f(z)|:\|f\|_{\mathfrak{F}_{\alpha,I}^p}\le 1\right\}\leq \displaystyle\sup_{z\in \mathbb{B}_I}\left\{|f_1(z)|:\|f_1\|_{\mathfrak{F}_{\alpha,\mathbb{C}_I}^p}\le 1\right\}+ \displaystyle\sup_{z\in \mathbb{B}_I}\left\{|f_2(z)|:\|f_2\|_{\mathfrak{F}_{\alpha,\mathbb{C}_I}^p}\le 1\right\}.$$ Let $q=x+yI_q$  and $z=x+yI.$  By using triangle inequality and   Theorem \ref{eq:115}, we have
  $$|f(q)|\leq |f(z)|+|f(\bar z)|.$$ On taking supremum over all $q\in \mathbb{B},$ we conclude  that
$$\begin{array}{ccl}
\displaystyle\sup_{q\in \mathbb{B}}\left\{|f(q)|:\|f\|_{\mathfrak{F}_{\alpha}^p}\leq 1\right\}&\leq& \displaystyle\sup_{z\in \mathbb{B}_I}\left\{|f(z)|:\|f\|_{\mathfrak{F}_{\alpha,I}^p}\le 1\right\} + \displaystyle\sup_{ z\in \mathbb{B}_I}\left\{|f(\bar z)|:\|f\|_{\mathfrak{F}_{\alpha,I}^p}\le 1\right\}\\
&=&\displaystyle2 e^{\frac{\alpha}{2}|z|^2}\\
&=&2\displaystyle e^{\frac{\alpha}{2}|q|^2}.
\end{array}$$

 Hence the result.
\endpf
%%%%%$$$$$$$$$$$$$$$$$$$$$$$$$$$$$$$$$$$$$$$$$$$$$$$$$$$$$$$$$$$$$$$$$$$$$$$$$$$$$$$$$$$$$$$$$$$$$$$$$$$$
%%%%%%%%%%###################################################################
\begin{corollary}
Suppose  $p>1$ and $\alpha >0.$ If   $f$  is in $\mathfrak{F}_{\alpha}^p( \mathbb{B}),$   then $$|f(q)|\leq  2^{p+1}e^{\frac{\alpha}{2}|q|^2}\|f\|_{\mathfrak{F}_{\alpha}^p}, for~all~q=x+yI_q\in \mathbb{B}.$$
\end{corollary}
\beginpf
  Let  $f\in \mathfrak{F}_{\alpha,I}^p.$   Then by Remark \ref{eq:100} and  \cite[Corollary 2.8]{zhu12}, we have  
\begin{eqnarray}\label{eq:152}
|f(z)|^p&\leq& 2^{p-1}(|f_1(z)|^p+|f_2(z)|^p) \nonumber\\
&\leq& 2^{p-1}( e^{\frac{\alpha p}{2}|z|^2}\|f_1\|^p_{\mathfrak{F}_{\alpha, \mathbb{C}_I}^p}+ e^{\frac{\alpha p}{2}|z|^2}\|f_2\|^p_{\mathfrak{F}_{\alpha, \mathbb{C}_I}^p})\nonumber\\
&\leq&2^p  e^{\frac{\alpha p}{2}|z|^2}\|f\|^p_{\mathfrak{F}_{\alpha, I}^p}.\\
 \nonumber\end{eqnarray}
Now, on applying Representation Formula,  condition (\ref{eq:152}) and  Remark \ref{eq:160}, we obtain $$|f(q)|^p\leq 2|f(z)|^p\leq    2^{p+1}  e^{\frac{\alpha p}{2}|z|^2}\|f\|^p_{\mathfrak{F}_{\alpha, I}^p}\leq 2^{p+1}  e^{\frac{\alpha p}{2}|z|^2}\|f\|^p_{\mathfrak{F}_{\alpha}^p.}$$
\endpf
%%%%%%$$$$$$$$$$$$$$$$$$$$$$$$$$$$$$$$$$$$$$$$$$$$$$$$$$$$$$$$$$$$$$$$$$$$$
%%%%%################################################################
We can easily prove the following result. 
\begin{proposition}
Let $1<p<\infty$ and $r\in (0,1).$ Then for any  $f\in \mathfrak{F}_{\alpha}^p$  $$\lim_{r\to 1}\|f_r-f\|^p_{\mathfrak{F}_{\alpha}^p}=0,$$ where $f_r(q)=f(rq)=\displaystyle\sum_{k=0}^\infty r^kq^ka_k,$ for all $q\in \mathbb{B}.$ 
\end{proposition}
\beginpf
Let $f\in \mathfrak{F}_{\alpha}^p.$  Then $f\in \mathfrak{F}_{\alpha,I}^p.$ Let $I,J \in \mathbb{S}$ be such that $I\bot J.$ Let $f_1, f_2$ be holomorphic functions in   $ \mathbb{B}_I.$ By Remark \ref{eq:100},  it follows that $f_1,f_2$ lie in the complex Fock space $\mathfrak{F}_{\alpha, \mathbb{C}_I}^p.$ By applying corresponding results   \cite[Proposition  2.9 (a)]{zhu12} to $f_1, f_2$ in  $\mathfrak{F}_{\alpha, \mathbb{C}_I}^p,$ we obatin    $\displaystyle{\lim_{r\to 1}}\|f_{l,r}-f_l\|^p_{\mathfrak{F}_{\alpha, \mathbb{C}_I}^p}=0,l=1,2.$  Since $ \|f\|^p_{\mathfrak{F}_{\alpha}^p}\leq 2^p \|f\|^p_{\mathfrak{F}_{\alpha,I}^p},$ we have
$$\begin{array}{ccl}
\displaystyle\lim_{r\to 1}\|f_r-f\| ^p_{\mathfrak{F}_{\alpha}^p}&\leq& 2^p\displaystyle\lim_{r\to 1}\|f_r-f\|^p_{\mathfrak{F}_{\alpha,I}^p}\\
&\leq&2^p\displaystyle\left(\lim_{r\to 1}\|f_{1,r}-f_1\|_{\mathfrak{F}_{\alpha, \mathbb{C}_I}^p}-\lim_{r\to 1}\|f_{2,r}-f_2\|_{\mathfrak{F}_{\alpha, \mathbb{C}_I}^p}\right)\\
&=&0.\\
\end{array}$$
\endpf
%%%%%%%$$$$$$$$$$$$$$$$$$$$$$$$$$$$$$$$$$$$$$$$$$$$$$$$$$$$$$$$$$$$$$$$$$$$$$$$$$
%%%%%%%%%%%%%@@@@@@@@@@@@@@@@@@@@@@@@@@@@@@@@@@@@@@@@@@@@@@@@@@@@
\begin{proposition}
For $1<p<\infty,$ the slice regular Fock space is the  closure of the sequence $\{p_m\}$ of quaternionic polynomials of the form $p_m(q)=\displaystyle\sum_{k=0}^m q^k\beta_{m,k},$ where $\beta_{m,k}\in \mathbb{H}$  with norm $\|.\|_{ \mathfrak{F}_{\alpha}^p}.$  In particular,  the slice regular Fock  space $\mathfrak{F}_{\alpha}^p$ is separable.
\end{proposition}
\beginpf
Suppose $f\in  \mathfrak{F}_{\alpha}^p.$ Then $f\in  \mathfrak{F}_{\alpha,I}^p$ so that $f_1, f_2 \in hol(\mathbb{B}_I),$ where $f_l, l=1,2,$ is given by Splitting Lemma \ref{eq:1}.  Let $\beta_{m,k}=\zeta_{m,k}+\gamma_{m,k}J,$ where  $\zeta_{m,k}, \gamma_{m,k}\in \mathbb{C}_I.$ By denseness property of polynomials in complex Fock space, we can choose  polynomials of the form $p_{1,m}(z)=\displaystyle\sum_{k=0}^m z^k\zeta_{m,k}~~ \mbox{and}~~  p_{2,m}(z)=\displaystyle\sum_{k=0}^m z^k\gamma_{m,k}.$ Applying  \cite[Proposition  2.9 (b)]{zhu12} to each $f_l, p_{l,m}, l=1,2,$ we see $\|p_{l,m}-f_l\|_{  \mathfrak{F}_{\alpha,\mathbb{C}_I}^p}\to 0$ as $m\to \infty.$ Thus, we have 
 $$\begin{array}{ccl}
\|f-p_m\|_{ \mathfrak{F}_{\alpha}^p}&\leq&2^p \|f-p_m\|_{ \mathfrak{F}_{\alpha,I}^p}\\
&=&2^p \|(f_1+f_2J)-(p_{1,m}+p_{2,m}J)\|_{ \mathfrak{F}_{\alpha, \mathbb{C}_I}^p}\\
&\leq&2^p \|f_1-p_{1,m}\|_{ \mathfrak{F}_{\alpha,\mathbb{C}_I}^p}-2^p \|f_2-p_{2,m}\|_{ \mathfrak{F}_{\alpha,\mathbb{C}_I}^p}\to 0~ \mbox{as}~ m\to \infty.
\end{array}$$
Hence,  $\mathfrak{F}_{\alpha}^p$  is separable.  
\endpf
%%%%%%$$$$$$$$$$$$$$$$$$$$$$$$$$$$$$$$$$$$$$$$$$$$$$$$$$$$$$$$$$$$$$$$$$$$$$$$$$$$$$$$
%%%%%%%$###################################################################
\begin{proposition}
For $1<p<u<\infty$  with $\frac{1}{p}+\frac{1}{u}=1, \; \mathfrak{F}_{\alpha}^p\subset  \mathfrak{F}_{\alpha}^u.$ Moreover $\|f\|_{ \mathfrak{F}_{\alpha}^u}^u\leq \displaystyle2^{u+1}\frac{u}{p} \|f\|_{ \mathfrak{F}_{\alpha}^p}^u.$
\end{proposition}
\beginpf
Let  $f\in  \mathfrak{F}_{\alpha}^p.$ For any $I,J \in \mathbb{S}$ with $I\bot J.$ Then Lemma  \ref{eq:1}, guarantees  the existence of holomorphic functions $f_1,f_2: \mathbb{B}\cap \mathbb{C}_I\to \mathbb{C}_I$ such that $Q_I[f](z)=f_1(z)+f_2(z)J,$ for all $z=x+yI\in \mathbb{B}_I. $  From Remark \ref{eq:100}, it follows that $f_1,f_2$ lie in the complex Fock space  $\mathfrak{F}_{\alpha, \mathbb{C}_I}^p.$ Therefore, from 
\cite[Theorem  2.10]{zhu12}, we have $\|f_l\|^u_{ \mathfrak{F}_{\alpha, \mathbb{C}_I}^u}\leq \frac{u}{p}  \|f_l\|^u_{ \mathfrak{F}_{\alpha, \mathbb{C}_I}^p}$ for $l=1,2.$ Furthermore,
\begin{eqnarray}\label{eq:154}
\displaystyle\frac{\alpha u}{2\pi}\displaystyle\int_{\mathbb{B}_I}\left|f(z)e^{\frac{-\alpha}{2}|z|^2}\right|^udA_I(z)&\leq&   \displaystyle 2^{u-1}\frac{\alpha u}{2\pi}\displaystyle\int_{\mathbb{B}_I}\left|f_1(z)e^{\frac{-\alpha}{2}|z|^2}\right|^udA_I(z) \nonumber\\
&+& \displaystyle 2^{u-1}\frac{\alpha u }{2\pi}\displaystyle\int_{\mathbb{B}_I}\left|f_2(z)e^{\frac{-\alpha}{2}|z|^2}\right|^udA_I(z) \nonumber\\
&=&\displaystyle 2^{u-1}(\|f_1\|_{ \mathfrak{F}_{\alpha, \mathbb{C}_I}^u}^u+ \|f_2\|_{ \mathfrak{F}_{\alpha, \mathbb{C}_I}^u}^u) \nonumber\\
&\leq&\displaystyle 2^{u-1}\frac{u}{p}(\|f_1\|_{ \mathfrak{F}_{\alpha, \mathbb{C}_I}^p}^u+ \|f_2\|_{ \mathfrak{F}_{\alpha, \mathbb{C}_I}^p}^u) \nonumber\\
&\leq&\displaystyle2^u\frac{u}{p}\|f\|_{ \mathfrak{F}_{\alpha,I}^p}^u.
\end{eqnarray}
Now, Theorem \ref{eq:115}  follows $|f(q)|\leq |f(z)|+ |f(\bar z)|,$ where $q=x+yI_q  \in \mathbb{B}$  with $I_q=\frac{Im(q)}{|Im(q)|}$ and $z=x+yI \in \mathbb{B}_I$ for all $x,y\in \mathbb{R}$  and by equation (\ref{eq:154}), we conclude that
$$\begin{array}{ccl}
\displaystyle\frac{\alpha u}{2\pi}\displaystyle\int_{\mathbb{B}_I}\left|f(q)e^{\frac{-\alpha}{2}|q|^2}\right|^udA_I(q)&\leq&   \displaystyle\frac{\alpha u}{2\pi}\displaystyle\int_{\mathbb{B}_I}\left|f(z)e^{\frac{-\alpha}{2}|z|^2}\right|^udA_I(z)\\
&+& \displaystyle\frac{\alpha u}{2\pi}\displaystyle\int_{\mathbb{B}_I}\left|f(\bar z)e^{\frac{-\alpha}{2}|z|^2}\right|^udA_I(\bar z)\\
&\leq& 2\displaystyle\frac{\alpha u}{2\pi}\displaystyle\int_{\mathbb{B}_I}\left|f(z)e^{\frac{-\alpha}{2}|z|^2}\right|^udA_I(z)\\
&\leq&\displaystyle \displaystyle2^{u+1}\frac{u}{p}\|f\|_{ \mathfrak{F}_{\alpha,I}^p}^u\\
&\leq&\displaystyle2^{u+1}\frac{u}{p}\|f\|_{ \mathfrak{F}_{\alpha}^p}^u.\\
\end{array}$$
Hence   the  result.
\endpf
%%%%%$$$$$$$$$$$$$$$$$$$$$$$$$$$$$$$$$$$$$$$$$$$$$$$$$$$$$$$$$$$$$$$$$$$$$$$$$$$$$
%%%%%%%%%^^^^^^^^^^^^^^^^^^^^^^^^^^^^^^^^^^^^^^^^^^^^^^^^^^^^^^^^^^^^^^^^^^
\begin{proposition}
Let $1<p<\infty.$ Then  for all $q,w \in \mathbb{B},$ the  function  $f(q)=\displaystyle P_I\sum_{m=0}^n e^{\beta q \bar w_m} a_m$ is dense in $ \mathfrak{F}_{\alpha}^p$ for some positive parameters $\alpha$ and $\beta.$
\end{proposition}
\beginpf
If $f\in \mathfrak{F}_{\alpha}^p,$ then $f\in \mathfrak{F}_{\alpha,I}^p.$ Let $f_1, f_2 \in hol(\mathbb{B}_I)$  given as in Lemma \ref{eq:1}  such that  $Q_I[f]=f_1+f_2J.$ Therefore from  \cite[Lemma 2.11]{zhu12},  each  functions of the form   $f_1(z)=\displaystyle\sum_{m=0}^n e^{\beta q \bar w_m} c_m$ and $f_2(z)=\displaystyle\sum_{m=0}^n e^{\beta q \bar w_m} d_m$ is dense on $ \mathfrak{F}_{\alpha, \mathbb{C}_I}^p$ on $\mathbb{B}_I.$ Consequently, $$Q_I[f](z)=f_1(z)+f_2(z)J=\sum_{m=0}^n e^{\beta z \bar w_m} c_m+\sum_{m=0}^n e^{\beta z \bar w_m} d_m J.$$ This implies that for each $q,w\in \mathbb{B},$ we have $$f=P_I \circ Q_I[f]=P_I\left[\sum_{m=0}^n e^{\beta q \bar w_m}(c_m+d_mJ)\right]=P_I\left[\sum_{m=0}^n e^{\beta q \bar w_m}a_m\right],$$ where the sequence $a_m=c_m+d_mJ$ lie in $l^p(\mathbb{H}).$ Thus, the  density  of  $f_1$ and $f_2$ in $\mathfrak{F}_{\alpha,\mathbb{C}_I}^p$ implies  $f$ is dense in $\mathfrak{F}_{\alpha,I}^p.$ Therefore, from Proposition \ref{eq:153}, we conclude that the set of functions $f$ is dense $\mathfrak{F}_{\alpha}^p.$
\endpf
%%%%%%%%%$$$$$$$$$$$$$$$$$$$$$$$$$$$$$$$$$$$$$$$$$$$$$$$
%%%%%%%%%%%#############################################################
\begin{proposition}
Let $0<p\leq \infty$ and  for some  $I\in \mathbb{S}. $ Then  $f\in \mathfrak{F}_{\alpha}^p({\mathbb{B}})$ if and only if there exists $\mathbb{H}$-valued Borel measure $\mu$ such that 
\beq \label{eq:12}
f(q)=\displaystyle\int_{\mathbb{B}_I}e^{\alpha \bar \zeta q -\frac{\alpha}{2}| \zeta|^2}d\mu( \zeta)\mbox{ for ~each} ~  \zeta,q \in \mathbb{B} ~  and ~\{|\mu|(G_1+w): w\in \mathbb{H}\}\in l^p(\mathbb{H}).
\eeq
\end{proposition}
\beginpf
Suppose $f\in \mathfrak{F}_{\alpha}^p({\mathbb{B}})$ implies $f\in \mathfrak{F}_{\alpha,I}^p({\mathbb{B}_I}).$ Let $J\in \mathbb{S}$  be such that $J\bot I.$ Then $f$ decomposes as $f_I=f_1+f_2J,$ where  $f_1,f_2:\mathbb{B}\cap \mathbb{C}_I \to \mathbb{C}_I$ with  $J\bot I.$  Clearly the holomorphic functions $f_1, f_2$ lie in the complex Fock space $\mathfrak{F}_{\alpha,\mathbb{C}_I}^p$ on $\mathbb{B}_I.$ Further, for each $f_l\in \mathfrak{F}_{\alpha,\mathbb{C}_I}^p(\mathbb{B}_I), l=1,2$ (see \cite[ p 91]{zhu12}), there exist complex  positive Borel measure $\mu_1$ and $\mu_2$ on $\mathbb{B}_I$ such that 
$f_l(z)=\displaystyle\int_{\mathbb{B}_I}e^{\alpha \bar  \zeta z -\frac{\alpha}{2}| \zeta|^2}d\mu_l( \zeta),~\mbox{ for each} ~ z=x+yI \in \mathbb{B}_I ~ \mbox{and} ~~\{|\mu_l|(G_r+w): w\in r\mathbb{R}^2\}\in l^p(\mathbb{C}_I).$ Now if we decompose   $\mu=\mu_1+\mu_2J,$  then, we can write
$$\begin{array}{ccl}
f(q)= Q_I[f_1+f_2](q)&=&\displaystyle\int_{\mathbb{B}_I}e^{\alpha \bar  \zeta q -\frac{\alpha}{2}| \zeta|^2}d\mu_1( \zeta)+\displaystyle\int_{\mathbb{B}_I}e^{\alpha \bar  \zeta q -\frac{\alpha}{2}| \zeta|^2}d\mu_2( \zeta)J\\
%&=&\displaystyle\int_{\mathbb{B}_i}e^{\alpha \bar  \zeta q -\frac{\alpha}{2}| \zeta|^2}(d\mu_1(a)+d\mu_2(a)j)\\
&=&\displaystyle\int_{\mathbb{B}_I}e^{\alpha \bar  \zeta q -\frac{\alpha}{2}| \zeta|^2}d\mu( \zeta).
\end{array}$$
Conversely, assume the condition (\ref{eq:12}) holds. So we can find complex valued Borel measure $\mu_1$ and $\mu_2$ in $\mathbb{C}_I$ such that $\mu=\mu_1+\mu_2J.$ Therefore, for each  $z\in \mathbb{B}_I$$$f_1(z)+f_2(z)J=Q_I[f](z)=\displaystyle\int_{\mathbb{B}_I}e^{\alpha \bar  \zeta z -\frac{\alpha}{2}| \zeta|^2}d\mu_1( \zeta)+\displaystyle\int_{\mathbb{B}_I}e^{\alpha \bar  \zeta z -\frac{\alpha}{2}| \zeta|^2}d\mu_2( \zeta)J.$$ Therefore,  $f_l(z)=\displaystyle\int_{\mathbb{B}_I}e^{\alpha \bar  \zeta  z -\frac{\alpha}{2}| \zeta|^2}d\mu_l( \zeta),l=1,2$ and as   $\{|\mu|(G_r+w): w\in \mathbb{H}\}\in l^p(\mathbb{H}),$ it follows that $\{|\mu_l|(G_r+w): w\in r\mathbb{R}^2\}\in l^p(\mathbb{C}_I).$ This implies $f_1$ and $f_2$ belong to complex Fock space $\mathfrak{F}_{\alpha,\mathbb{C}_I}^p(\mathbb{B}_I)$ which is equivalent to  $f \in \mathfrak{F}_{\alpha,I}^p$ and hence $f\in \mathfrak{F}_{\alpha}^p$ in $\mathbb{B}.$
 \endpf

\end{document}